\documentclass[a4paper, 11pt] {article}
\usepackage[cp1251]{inputenc}
\usepackage[russian, english]{babel}
\usepackage{amsfonts}
\usepackage{mathtext}
\usepackage{cite}

\usepackage{graphicx}

\usepackage{verbatim}
\usepackage{amsmath}
\usepackage{amsthm}
\usepackage{amssymb}
\usepackage{delarray}
\usepackage{mathrsfs}
\usepackage{xcolor}

\begin{document}
%\vspace*{5mm}

\thispagestyle{empty}

\begin{center}
{\bf\large ON AN OPEN QUESTION IN RECOVERING STURM--LIOUVILLE-TYPE OPERATORS WITH DELAY}
\end{center}

\begin{center}
{\bf\large Neboj\v{s}a Djuri\'c\footnote{Faculty of Architecture, Civil Engineering and Geodesy, University
of Banja Luka, {\it nebojsa.djuric@aggf.unibl.org}} and Sergey Buterin\footnote{Department of Mathematics,
Saratov State University, {\it buterinsa@info.sgu.ru}}}
\end{center}

{\bf Abstract.} In recent years, there appeared a considerable interest in the inverse spectral theory for
functional-differential operators with constant delay. In particular, it is well known that specification of
the spectra of two operators $\ell_j,$ $j=0,1,$ generated by one and the same functional-differential
expression $-y''(x)+q(x)y(x-a)$ under the boundary conditions $y(0)=y^{(j)}(\pi)=0$ uniquely determines the
complex-valued square-integrable potential $q(x)$ vanishing on $(0,a)$ as soon as $a\in[\pi/2,\pi).$ For many
years, it has been a challenging {\it open question} whether this uniqueness result would remain true also
when $a\in(0,\pi/2).$ Recently, a positive answer was obtained for the case $a\in[2\pi/5,\pi/2).$ In this
paper, we give, however, a {\it negative} answer to this question for $a\in[\pi/3,2\pi/5)$ by constructing an
infinite family of iso-bispectral~potentials. Some discussion on a possibility of constructing a similar
counterexample for other types of boundary conditions is provided, and new open questions are outlined.


\\

{\large\bf 1. Introduction}
\\

One of the first results in the inverse spectral theory says that the spectra of two boundary value problems
for one and the same Sturm--Liouville equation with one common boundary condition:
\begin{equation}\label{0}
-y''(x)+q(x)y(x)=\lambda y(x), \quad y(0)=y^{(j)}(\pi)=0, \quad j=0,1,
\end{equation}
uniquely determine the potential $q(x),$ see \cite{B}, where also local solvability and, actually, stability
of this inverse problem were established in the class of real-valued $q(x)\in L_2(0,\pi).$ Later on, these
results were refined and generalized to other classes of potentials and boundary conditions \cite{Kar, MO,
FY01, HM04, SavShk, BK19}. Moreover, there appeared methods, which gave global solution for the inverse
Sturm--Liouville problem as well as for inverse problems for other classes of differential operators (see,
e.g., monographs \cite{FY01, M, L, Yur02}).

In recent years, there was a considerable interest in inverse problems also for Sturm--Liouville-type
operators with deviating argument (see \cite{Pik91, FrYur12, VladPik16, ButYur19, ButPikYur17, Ign18,
BondYur18-1, BondYur18-2, VPV19, Yur19, DV19, Yur19-2, Yur20, SS19, Dur, Yang, Wang, BK, Wang20} and
references therein), which are often more adequate for modelling various real world processes frequently
possessing a nonlocal nature.

For $j=0,1,$ denote by $\{\lambda_{n,j}\}$ the spectrum of the boundary value problem ${\cal L}_j(a,q)$ of
the form
\begin{equation}\label{1}
-y''(x)+q(x)y(x-a)=\lambda y(x),\quad 0<x<\pi,
\end{equation}
\begin{equation}\label{2}
y(0)=y^{(j)}(\pi)=0,
\end{equation}
where $q(x)$ is a complex-valued function in $L_2(0,\pi)$ and $q(x)=0$ on $(0,a),$ while $a\in(0,\pi).$ For
functional-differential operators as well as for other classes of nonlocal ones, classical method of the
inverse spectral theory for differential operators do not work. Consider the following inverse problem.

\medskip
{\bf Inverse Problem 1.} Given the spectra $\{\lambda_{n,0}\}$ and $\{\lambda_{n,1}\},$ find the potential
$q(x).$

\medskip
Instead of (\ref{2}), one could alternatively impose the boundary conditions
\begin{equation}\label{rob}
y'(0)-hy(0)=y^{(j)}(\pi)=0
\end{equation}
with a complex parameter $h,$ which can be found given the input data of Inverse Problem~1. Various aspects
of Inverse Problem~1 were studied in \cite{Pik91, FrYur12, VladPik16, ButYur19, ButPikYur17, Ign18,
BondYur18-1, BondYur18-2, VPV19, Yur19, DV19, Yur19-2, Yur20, Dur} and other works.

The values $\pi/2$ and $2\pi/5$ of the delay parameter $a$ are critical for the solution of this inverse
problem. Specifically, in the case $a\ge\pi/2,$ the dependence of the characteristic function of each problem
${\cal L}_j(a,q)$ on the potential $q(x)$ is linear, while for $a<\pi/2$ it is nonlinear. In the linear case,
the uniqueness of solution of Inverse Problem~1 is well known. Moreover, this problem is then overdetermined.
Thus, in \cite{ButYur19}, the conditions on an arbitrary increasing sequence of natural numbers
$\{n_k\}_{k\ge1}$ were obtained that are necessary and sufficient for the unique determination of the
potential $q(x)$ by specifying the corresponding subspectra $\{\lambda_{n_k,0}\}$ and $\{\lambda_{n_k,1}\}.$
The nonlinear case $a<\pi/2$ is much more difficult. However, in \cite{Ign18, FrYur12, ButPikYur17} for
various types of boundary conditions and for any fixed $a\in(0,\pi),$ it was established that if the given
spectra coincide with the spectra of the corresponding problems with the zero potential, then $q(x)$ is zero
too.

Under the general settings, i.e. for arbitrary nonzero potentials, it was a long term challenging {\it open
question} whether the uniqueness of solution of Inverse Problem~1 takes place for $a<\pi/2.$ Recently, a
positive answer to this question was given independently in \cite{BondYur18-1} (for the conditions (\ref{2}))
and in \cite{VPV19} (for (\ref{rob})) as soon as $a\in[2\pi/5,\pi/2).$

In the present paper, we give a {\it negative} answer to the {\it open question} formulated above in the case
when $a\in[\pi/3,2\pi/5).$ Namely, for each such $a,$ we construct an infinite family of different
iso-bispectral potentials~$q(x),$ i.e. of those for which the problems ${\cal L}_0(a,q)$ and ${\cal
L}_1(a,q)$ have one and the same pair of spectra. This looks quite unusual in light of the classical results
for the case of $a=0.$ Our counterexample appears even more unexpected under consideration that the recent
paper \cite{Yur20} announces that specification of the spectra of both boundary value problems consisting
of~(\ref{1}) and (\ref{rob}) uniquely determines the potential $q(x)$ also for $a\in[\pi/3,2\pi/5).$ Although
in \cite{Yur20} Robin boundary conditions were imposed also in the point $\pi,$ they can be easily reduced
to~(\ref{rob}). In Section~3, we discuss why our result does not refute the uniqueness theorem in
\cite{Yur20} for $a\in[\pi/3,2\pi/5)$ (see Remark~2).

In the next section, we formulate the main result of the paper (Theorem~1) and outline some new prospects in
studying Inverse Problem~1 (Remark~1). The proof of Theorem~1 is provided in Section~3.
\\

{\large\bf 2. The main result}
\\

Fix $a\in[\pi/3,2\pi/5).$ For a nonzero real-valued function $h(x)\in L_2(5a/2,\pi),$ we consider the
integral operator
\begin{equation}\label{2.1}
M_hf(x)=\int_\frac{3a}2^{\pi-x+\frac{a}2} K_h\Big(x+t-\frac{a}2\Big)f(t)\,dt, \;\; \frac{3a}2<x<\pi-a, \quad
{\rm where} \quad K_h(x)=\int_x^\pi h(\tau)\,d\tau.
\end{equation}
Thus, $M_h$ is a nonzero compact Hermitian operator acting in the space $L_2(3a/2,\pi-a).$ Hence, it has at
least one nonzero real eigenvalue $\eta.$ We assume that $\eta=1,$ which can always be achieved by choosing
$h(x)/\eta$ instead of $h(x).$ Let $e(x)$ be the corresponding eigenfunction,~i.e.
\begin{equation}\label{2.2}
M_he(x)=e(x), \quad \frac{3a}2<x<\pi-a.
\end{equation}
Consider the one-parametric family of potentials $B:=\{q_\alpha(x)\}_{\alpha\in{\mathbb C}}$ determined by
the formula
\begin{equation}\label{2.3}
q_\alpha(x)=\left\{\begin{array}{cl}\displaystyle 0,
 &\displaystyle x\in\Big(0,\frac{3a}2\Big)\cup(\pi-a,2a)\cup\Big(\pi-\frac{a}2,\frac{5a}2\Big),\\[3mm]
\displaystyle \alpha e(x), &\displaystyle x\in\Big(\frac{3a}2,\pi-a\Big),\\[3mm]
%\displaystyle 0, &\displaystyle x\in(\pi-a,2a),\\[3mm]
\displaystyle -\alpha K_h\Big(x+\frac{a}2\Big)\int_\frac{3a}2^{x-\frac{a}2} e(t)\,dt,
 &\displaystyle x\in\Big(2a,\pi-\frac{a}2\Big),\\[3mm]
%\displaystyle 0, &\displaystyle x\in\Big(\pi-\frac{a}2,\frac{5a}2\Big),\\[3mm]
h(x), &\displaystyle x\in\Big(\frac{5a}2,\pi\Big).
\end{array}\right.
\end{equation}

The main result of the present paper is the following theorem.

\medskip
{\bf Theorem 1. }{\it For $j=0,1,$ the spectrum $\{\lambda_{n,j}\}$ of the boundary value problem ${\cal
L}_j(a,q_\alpha)$ does not depend on $\alpha.$}

\medskip
Theorem~1 means that the problems ${\cal L}_0(a,q_\alpha)$ and ${\cal L}_1(a,q_\alpha)$ have one and the same
pair of spectra $\{\lambda_{n,0}\}$ and $\{\lambda_{n,1}\}$ for all values of the parameter
$\alpha\in{\mathbb C},$ i.e. the solution of Inverse Problem~1 is not unique. We note that, since the zero
potential does not belong to $B,$ Theorem~1 does not contradict the uniqueness result in \cite{FrYur12}. The
proof of Theorem~1 is given in the next section.

\medskip
{\bf Remark 1.} This reminds the situation of the following operator with frozen argument:
$$
\ell y=-y''(x)+q(x)y(b), \quad y^{(\alpha)}(0)=y^{(\beta)}(1)=0, \quad \alpha,\beta\in\{0,1\},\quad
b\in[0,1].
$$
It is known that the unique recoverability of $q(x)$ from the spectrum of $\ell$ depends on the triple of
parameters $(\alpha,\beta,b).$ In \cite{BK}, one can find a complete description of all degenerate and
non-degenerate cases for rational values of $b,$ while for irrational ones the uniqueness always takes place,
see \cite{Wang20}.

Theorem~1 opens the analogous type of questions for Inverse Problem~1, which consist in giving description of
ranges of the delay parameter $a$ along with the types of boundary conditions, for which the uniqueness of
solution of Inverse Problem~1 takes place. In the opposite degenerate cases, it would be interesting to
describe classes of iso-bispectral potentials. It is especially important to investigate the case of
arbitrarily small values of $a$ making the problem ''close'' to the classical case~(\ref{0}).
\\

{\large\bf 3. Proof of Theorem~1}
\\

Let $S(x,\lambda)$ be a solution of equation (\ref{1}) under the initial conditions $S(0,\lambda)=0$ and
$S'(0,\lambda)=1.$ For $j=0,1,$ eigenvalues of the problem ${\cal L}_j(a,q)$ coincide with zeros of its
characteristic function $\Delta_j(\lambda)=S^{(j)}(\pi,\lambda).$ Thus, the spectrum of any problem ${\cal
L}_j(a,q)$ does not depend on $q(x)\in B$ for some subset $B\subset L_2(0,\pi)$ as soon as neither does the
corresponding characteristic function $\Delta_j(\lambda).$ Put $\rho^2=\lambda$ and denote
\begin{equation}\label{3.0}
\omega:=\int_a^\pi q(x)\,dx.
\end{equation}
The following representations hold (see, e.g., \cite{Dur}):
\begin{equation}\label{3.1}
\Delta_0(\lambda)=\frac{\sin\rho\pi}\rho -\omega\frac{\cos\rho(\pi-a)}{2\rho^2} +\frac1{2\rho^2}\int_a^\pi
w_0(x)\cos\rho(\pi-2x+a)\,dx,
\end{equation}
\begin{equation}\label{3.2}
\Delta_1(\lambda)=\cos\rho\pi +\omega\frac{\sin\rho(\pi-a)}{2\rho} -\frac1{2\rho}\int_a^\pi
w_0(x)\sin\rho(\pi-2x+a)\,dx,
\end{equation}
where the function $w_0(x)$ is determined by the following formula for $k=0:$
\begin{equation}\label{3.5}
w_k(x)=\left\{\begin{array}{l}
\displaystyle q(x),\quad x\in\Big(a,\frac{3a}2\Big)\cup\Big(\pi-\frac{a}2,\pi\Big),\\[3mm]
\displaystyle q(x)+Q_k(x),\quad x\in\Big(\frac{3a}2,\pi-\frac{a}2\Big),
\end{array}\right.
\end{equation}
while
\begin{equation}\label{3.6}
Q_k(x)=\int_{x-\frac{a}2}^{\pi-a} \Big(q(t+a)\int_a^t q(\tau)\,d\tau -q(t)\int_{t+a}^\pi q(\tau)\,d\tau
-(-1)^k\int_{t+a}^\pi q(\tau-t)q(\tau)\,d\tau\Big)dt.
\end{equation}

Note that, since the function $\Delta_0(\lambda)$ is entire in $\lambda,$ representation (\ref{3.1}) implies
\begin{equation}\label{3.6.1}
\omega=\int_a^\pi w_0(x)\,dx,
\end{equation}
which can also be checked independently by direct calculation using (\ref{3.5}) and (\ref{3.6}) for $k=0.$
Thus, the spectrum of ${\cal L}_j(a,q),\,j=0,1,$ is independent of $q(x)\in B$ if so is the
function~$w_0(x).$

By changing the order of integration in (\ref{3.6}), it is easy to obtain the representation
\begin{equation}\label{3.7}
Q_k(x)=\int_a^{x-\frac{a}2} q(t)dt\int_{x+\frac{a}2}^\pi q(\tau)\,d\tau -(-1)^k\int_a^{\pi-x+\frac{a}2}
q(t)\,dt\int_{x+t-\frac{a}2}^\pi q(\tau)\,d\tau,
\end{equation}
where $x\in(3a/2,\pi-a/2).$ Let $q(x)=0$ on $(a,3a/2).$ Then (\ref{3.7}) takes the form
\begin{equation}\label{3.8}
Q_k(x)=\left\{\begin{array}{cl} \displaystyle -(-1)^k M_{q_2}q_1(x),
& \displaystyle x\in\Big(\frac{3a}2,\pi-a\Big),\\[3mm]
\displaystyle 0, & x\in(\pi-a,2a),\\[3mm]
\displaystyle K_{q_2}\Big(x+\frac{a}2\Big)\int_\frac{3a}2^{x-\frac{a}2} q_1(t)\,dt,
 & \displaystyle x\in\Big(2a,\pi-\frac{a}2\Big),
\end{array}\right.
\end{equation}
where
$$
q_1(x)=q(x), \;\; x\in\Big(\frac{3a}2,\pi-a\Big), \quad q_2(x)=q(x), \;\;x\in\Big(\frac{5a}2,\pi\Big),
$$
while $M_h$ and $K_h(x)$ are determined by (\ref{2.1}). Thus formulae (\ref{3.5}) and (\ref{3.8}) give
\begin{equation}\label{3.9}
w_k(x)=\left\{\begin{array}{cl} \displaystyle 0,
& \displaystyle x\in\Big(a,\frac{3a}2\Big),\\[3mm]
 \displaystyle q_1(x)-(-1)^k M_{q_2}q_1(x),
& \displaystyle x\in\Big(\frac{3a}2,\pi-a\Big),\\[3mm]
\displaystyle q(x), & x\in(\pi-a,2a),\\[3mm]
\displaystyle q(x)+K_{q_2}\Big(x+\frac{a}2\Big)\int_\frac{3a}2^{x-\frac{a}2} q_1(t)\,dt,
 & \displaystyle x\in\Big(2a,\pi-\frac{a}2\Big),\\[3mm]
\displaystyle q(x),
 & \displaystyle x\in\Big(\pi-\frac{a}2,\frac{5a}2\Big),\\[3mm]
 \displaystyle q_2(x),
 & \displaystyle x\in\Big(\frac{5a}2,\pi\Big).
\end{array}\right.
\end{equation}
Substituting $q(x)=q_\alpha(x)$ into (\ref{3.9}) for $k=0,$ where $q_\alpha(x)$ is determined by (\ref{2.3}),
and taking (\ref{2.2}) into account, we arrive at
$$
w_0(x)=\left\{\begin{array}{cl} \displaystyle 0, & \displaystyle x\in\Big(a,\frac{5a}2\Big),\\[3mm]
h(x), & \displaystyle x\in\Big(\frac{5a}2,\pi\Big).
\end{array}\right.
$$
Thus, according to (\ref{3.1}), (\ref{3.2}) and (\ref{3.6.1}), the characteristic function
$\Delta_j(\lambda)$ of the problem ${\cal L}_j(a,q_\alpha)$ for each $j=0,1$ is independent of $\alpha,$
which finishes the proof of Theorem~1. $\hfill\Box$

\medskip {%\color{blue}
{\bf Remark 2.} Thus, we have constructed a class $B$ of potentials on those the function $w_0(x)$ appearing
in representations (\ref{3.1}) and (\ref{3.2}) does not depend. By virtue of relation (\ref{3.6.1}), this
independence is inherited by both characteristic functions $\Delta_0(\lambda)$ and $\Delta_1(\lambda),$ which
gives Theorem~1.

Let us show why this strategy fails (at least, in the present form) in the case of boundary
conditions~(\ref{rob}). For $j=0,1,$ denote by ${\cal M}_j(a,h,q)$ the boundary value problem for equation
(\ref{1}) with boundary conditions~(\ref{rob}). Then eigenvalues of the problem ${\cal M}_j(a,0,q)$ coincide
with zeros of its characteristic function $\Theta_j(\lambda):=C^{(j)}(\pi,\lambda),$ where $C(x,\lambda)$ is
the solution of~(\ref{1}) under the initial conditions $C(0,\lambda)=1$ and $C'(0,\lambda)=0.$ Analogously to
(\ref{3.1}) and (\ref{3.2}), one can obtain the representations (see, e.g., \cite{VPV19}):
\begin{equation}\label{3.3}
\Theta_0(\lambda)=\cos\rho\pi +\omega\frac{\sin\rho(\pi-a)}{2\rho} +\frac1{2\rho}\int_a^\pi
w_1(x)\sin\rho(\pi-2x+a)\,dx,
\end{equation}
\begin{equation}\label{3.4}
\Theta_1(\lambda)=-\rho\sin\rho\pi +\frac\omega2 \cos\rho(\pi-a) +\frac12\int_a^\pi
w_1(x)\cos\rho(\pi-2x+a)\,dx,
\end{equation}
where $\omega$ is determined by formula (\ref{3.0}), while the function $w_1(x)$ has the form (\ref{3.5}) for
$k=1.$ Analogously to $B,$ one can construct a family $B_1$ of potentials $p_\alpha(x)$ on those the function
$w_1(x)$ does not depend. Indeed, for this purpose, the same scheme can be used but involving the operator
$-M_h$ instead of $M_h.$ However, the main difference from the case of boundary conditions (\ref{2}) is that,
in the case of (\ref{rob}), there is no relation analogous to (\ref{3.6.1}). In other words, the constant
$\omega$ is not determined by the function $w_1(x).$ Thus, each characteristic function (\ref{3.3}) and
(\ref{3.4}) could depend on $\alpha$ even if $w_1(x)$ did not. So the presented scheme does not refute the
uniqueness theorem in~\cite{Yur20}.
\\

{\bf Acknowledgement.} The first author was supported by the Project 19.032/961-103/19 of the Republic of
Srpska Ministry for Scientific and Technological Development, Higher Education and Information Society. The
second author was supported by Grants 19-01-00102 and 20-31-70005 of the Russian Foundation for Basic
Research.

\end{document}